\documentclass[11pt]{article}
\usepackage[all]{xy}
\usepackage{amstex}
\usepackage{delarray}
\usepackage{amssymb}
\newtheorem{defi}{ Definition :}[section]

\newtheorem{theo}[defi]{Theorem }
\newtheorem{lem}[defi]{Lemma }

\newtheorem{rem}[defi]{Remark }
\newcommand{\D}{\displaystyle}

\newcommand{\N}{{\mathbb N}}
\newcommand{\Z}{{\mathbb Z}}
\newcommand{\K}{{\mathbb K}}
\newcommand{\R}{{\mathbb R}}
\newcommand{\C}{{\mathbb C}}
\newcommand{\biv}{{\frac{\partial }{\partial x}
\wedge\frac{\partial }{\partial y}}}
\newcommand{\vx}{{\frac{\partial }{\partial x}}}
\newcommand{\vy}{{\frac{\partial }{\partial y}}}
\newcommand{\qh}{{quasihomogeneous }}
\newcommand{\degp}{{{\rm {d}} -\omega_1-\omega_2}}
\newcommand{\fct}{{{\cal F}(\K^2)}}
\newcommand{\chp}{{{\cal X}(\K^2)}}
\newcommand{\pois}{{{\cal V}(\K^2)}}
\def\qed{{\quad {\mbox {\tiny $\blacksquare$}}}}
\def\colon{{:}\;}
\def\|{|\;}

\def\endrem{}

\begin{document}
\addtocounter{footnote}{-1}

\begin{center}
{\bf {\Large Poisson cohomology in dimension two}}
\end{center}
\begin{center}
Philippe MONNIER
\end{center}
\vspace{0.5 cm}
\begin{abstract}
It is known that the computation of the Poisson cohomology is closely related to the 
classification of singularities of Poisson structures. 
In this paper, we will first 
look for the normal forms of germs at (0,0) of Poisson structures on 
$\K^2 \, (\K=\R {\mbox { or }} \C)$ and recall a result given by Arnold.
Then we will compute locally the Poisson cohomology of a particular type of Poisson
structure.
\end{abstract}

\section{Introduction}
The Poisson cohomology of a Poisson manifold gives several informations on the geometry of
the manifold. 
It was first introduced by Lichnerowicz in [L]. 
Unfortunately, the
computation of these cohomology spaces is quite complicated and few explicit results have 
been found.\\
In the symplectic case, Poisson cohomology is naturally isomorphic to de Rham 
cohomology.
The case of regular Poisson manifolds is discussed, for instance, in [V] and 
[X]. 
One can find some results on the Poisson cohomology of Poisson-Lie groups  
in [GW]. 
Some explicit computations are also done, 
for instance, in [Co1], [Co2] or [G].\\
In [Cr], Crainic links Poisson cohomology with the Morita equivalence.
Finally, one can find some discussion on Poisson cohomology and
Poisson homology in [ELW], [B] or [FT].\\
In the two-dimensional situation, some special cases on $\R^2$
have been studied. 
In [V], Vaisman began to compute the 
cohomology of $(x^2+y^2)\biv$. 
His idea was to consider the homomorphism 
$\iota^*\colon H^\bullet(\R^2)\to H^\bullet(\R^2 \setminus
\{(0,0)\})$ 
induced by the 
inclusion $\iota : \R^2 \setminus \{(0,0)\} \hookrightarrow \R^2$. 
A few years later, 
Nakanishi used this idea and computed the Poisson cohomology of quadratic Poisson structures
on $\R^2$ (see [N]).\\
In the present paper, our approach is more direct and uses some tools arising from the 
theory of singularities. 
More precisely, we first study (in section 3) the normal forms 
of the ``most interesting'' germs at (0,0) of Poisson structures vanishing at (0,0), 
and we rediscover the list given by Arnold in [A]. 
These normal forms are of type 
$$
\Pi=f(1+h)\biv\,,
$$ 
where $f$ and $h$ are \qh polynomials (there is a relation between 
their degrees).\\
Then, in section 4, we then compute locally the Poisson cohomology of Poisson structures 
of this type.\\
A generalisation of these results to the $n$-vectors on an 
$n$-dimensional manifold can be found in [Mo2].

\section{Quasi-homogeneity}
Throughout this text, $\K$ will indicate the field $\R$ or $\C$.\\
Let $(\omega_1,\omega_2)\in\N^*\times\N^*$.
We denote by $W$ the vector field $\omega_1x\vx +\omega_2y\vy$ on $\K^2$.
Now, let $T$ be a non zero p-vector (${\rm p}\in\{0,1,2\}$). 
We will say that $T$ is 
{\bf \qh} with weights $\omega_1,\omega_2$ and of (quasi)degree 
${\rm {d}}\in\Z$ if 
$$
{\cal L}_WT={\rm d}T\,,
$$ 
where ${\cal L}_W$ indicates the Lie 
derivative with respect to $W$. 
This condition can be written $[W,T]={\rm d}T$ 
where [.,.] indicates the Schouten bracket.
Note that $T$ is then polynomial.\\
If $f$ is a \qh polynomial of degree d then ${\rm d}=i\omega_1+j\omega_2$ with 
$(i,j)\in\N^2$; therefore, an integer is not necessarily the quasidegree of a polynomial.
If $f\in\K\big[[x,y]\big]$, we can write $f=\sum_{i=0}^{\infty}f_i$ with 
$f_i$ \qh  of degree $i$ (we adopt the convention that $f_i=0$ if $i$ is not a
quasidegree); $f$ is said to be of order d (${\mbox {ord}}(f)={\rm d}$) if all of its 
monomials have degree d or higher. 
For more details consult [AGV].\\
It is important to notice that $\vx$ is a  \qh vector field of degree $-\omega_1$
(in the same way $\deg (\vy)=-\omega_2$);
the minimal  degree of a vector field is  $-\max (\omega_1,\omega_2)$.
Note also that an integer can be the quasidegree of a vector field without being the 
quasidegree of a polynomial.
Finally, note that $\biv$ is \qh of degree $-\omega_1-\omega_2$.

\section{Local models of Poisson structures in dimension~2}
In the reference [A], Arnold gives a list of normal forms for Poisson structures on a
neighbourhood of $(0,0)$ in $\K^2$. 
In this section, we recall Arnold's theorem and we
give the idea of a proof which is similar to Arnold's 
(the approach is a little bit different).
For more details on this proof, consult [Mo1].\\
The particularity of the dimension two is that any 2-vector on a 2-dimensional
manifold is a Poisson structure. 
For an introduction to Poisson structures,
consult [CW] or [V].\\
The problem is the following: given $\Pi=F\biv$, a germ at 0 of Poisson structures on $\K^2$, 
we want to simplify its expression via a suitable local change of coordinates.\\

\noindent {\bf Notations :}
We denote by ${\cal F}(\K^2)$ \big(resp. ${\cal X}(\K^2), {\cal V}(\K^2)$\big) the vector
space of germs at (0,0) of (holomorphic if $\K=\C$ , analytic or ${\cal C}^\infty$
if $\K=\R$) functions (resp. vector fields, 2-vectors). 
We also denote
by ${Dif\!f}_0(\K^2)$
the group of local diffeomorphisms at (0,0) sending (0,0) to itself.
Finally, ${\cal F}_t(\K^2)$ \big(${\cal X}_t(\K^2), {\cal V}_t(\K^2)$\big) 
indicates the space of germs depending differentiably on $t\in\R$.\\

Two germs $\Pi=f\biv$ and $\Lambda=g\biv$ are called {\bf equivalent} if there 
exists $\varphi\in {Dif\!f}_0(\K^2)$ satisfying
$\varphi_*\Pi=\Lambda$. 
This condition yields
$g\circ\varphi=(Jac\, \varphi)f$ where $Jac\, \varphi$ indicates the
Jacobian of $\varphi$.\\
Two germs $f$ and $g$ are said to be {\bf R-equivalent} if there exists 
$\varphi\in {Dif\!f}_0(\K^2)$ satisfying $g\circ\varphi=f$.

Actually, it is not possible to give normal forms for every Poisson structures. 
We only
study Poisson structures determined by the germs of functions $F$ whose R-orbit is
``interesting'' enough. 
We will speak about it later.\\

The splitting theorem ([W]) allows us to assume that $\Pi_{(0,0)}=0$. 
Moreover, it is
quite easy to show that, if $F$ is regular at 0, then $\Pi$ is, up to a change of 
coordinates, the germ $x\biv$.

\begin{rem}
{\rm It is important to note that if two germs $f$ and $g$ are R-equivalent, then
the  germ $\Pi=f\biv$ will be equivalent to the germ of a Poisson
structure of type $\, ga\biv$ where $a(0,0)\neq 0$.}
\end{rem}
\endrem

Now, we consider germs of Poisson structures of type 
$$
\Pi=fa\biv\,,
$$ 
where $f$ vanishes and is singular at $(0,0)$, and $a(0,0)\neq 0$.\\
Moreover, we suppose that $f$ is a \qh polynomial of degree ${\rm d}>0$ with respect to 
$W=\omega_1x\vx+\omega_2y\vy$ ($\omega_1$ and $\omega_2$ are positive
integers). 
This additional assumption will be justified later.\\
Here and throughout, the quasihomogeneity will be with respect to $W$.\\
Arnold's theorem is the following:
\begin{theo} {\rm {\bf [A]}}
Up to a multiplicative constant, $\Pi$ is equivalent to the germ of a Poisson structure
of type $f(1+h)\biv$ where $h$ is a \qh polynomial of degree $\degp$ 
(on condition that $\degp$ is a quasidegree, otherwise this term disappears).
\label{fn poisson}
\end{theo}

It is possible to show (see [Mo1]) that $\Pi$ is equivalent to a germ of
Poisson structures of type $f(1+h+R)\biv$ where 
${\rm {ord}}\big(j_0^\infty(R)\big)>\degp$ ($j_0^\infty(R)$ indicates the $\infty$-jet of R
at (0,0)) and $h$ is a \qh polynomial of degree $\degp$.\\
Thus, in order to prove the former theorem, we have to ``remove'' the
term $R$. 
We are going
to use {\underline {Moser's path method}}.
For $t\in\R$, we put $\Pi_t=f(1+h+tR)\biv$ and we try to prove the existence of
$X_t\in {\cal X}_t(\K^2)$ satisfying 
$[X_t,\Pi_t]=-\frac{{\rm d}\Pi_t}{{\rm d} t}$.\\
Actually we will look for an $X_t$ of type $\alpha_t W$ with 
$\alpha_t\in{\cal F}_t(\K^2)$. \\
Then, if we put 
$R_t=\frac{R}{1+h+tR}$ and $\lambda_t=\degp +\frac{W.(h+tR)}{1+h+tR}$,
it is sufficient to prove the existence of $\alpha_t$ in ${\cal F}_t(\K^2)$, 
such that 
$$
W.\alpha_t -\lambda_t\alpha_t=R_t\qquad (E)\,.
$$
Let us note two things :\\
$\bullet$ if $\Pi$ is analytic (${\cal C}^\infty$) then so are $R_t$ and $\lambda_t$\\ 
$\bullet$ if $\degp$ is a quasidegree, then ${\rm {ord}}\big(j_0^\infty(R_t)\big)>\degp$.\\
Now, we just have to show that there exists $\alpha_t$ satisfying $(E)$.\\ 

\noindent {\bf Resolubility of equation $(E)$ :} 
The results we give here will be useful in the computation of the Poisson cohomology. 
That is why we are going to give quite detailed proofs of them.\\
We can write $\lambda_t=(\degp)+\mu_t$ 
where $\mu_t\in{\cal F}_t(\K^2)$ satisfies $\mu_t(0,0)=0$.
In order to show that $(E)$ admits a solution :\\
{\bf 1-} we prove that there exists $\beta_t\in{\cal F}_t(\K^2)$ satisfying
$W.\beta_t-\mu_t\beta_t=0$ with $\beta_t(0,0)\neq 0$\\
{\bf 2-} we prove that there exists $\gamma_t\in{\cal F}_t(\K^2)$ satisfying
$W.\gamma_t-(\degp)\gamma_t=\frac{R_t}{\beta_t}$\\
{\bf 3-} $\alpha_t:=\beta_t\gamma_t$ will then be a solution of $(E)$.\\ 

{\bf 1- } In order to show the first claim, we need the following result
whose proof can be found in [R].
\begin{theo}
Let $X_t$ be an element in ${\cal X}_t(\K^2)$ having an isolated
singularity at (0,0). 
Moreover, 
suppose that the eigenvalues of its linear component at (0,0) do not
vanish. 
Take
$h_t$ in ${\cal C}_t^\infty(\R^2)$ flat at (0,0). 
Then there exists 
$g_t\in{\cal C}_t^\infty(\R^2)$ flat at (0,0) satisfying $X_t.g_t=h_t$ for any $t$.
\label{roussarie}
\end{theo}
We also need the following lemma.
\begin{lem}
If $T_t\in{\cal F}_t(\K^2)$ satisfies $T_t(0,0)=0$, then there exists 
$\nu_t\in{\cal F}_t(\K^2)$  such that $W.\nu_t=T_t$.
\label{E1}
\end{lem}
{\sf Proof of the lemma }: \\
{\underline {Formal case :}} Assume that $T_t\in\K_t\big[[x,y]\big]$; we have 
$T_t=\sum_{i>0}T_t^{(i)}$ where $T_t^{(i)}$ is \qh of degree $i$. 
If we put
$\nu_t=\sum_{i>0}\frac{T_t^{(i)}}{i}$ we get $W.\nu_t=T_t$.\\
{\underline {Analytical case :}} Assume that $T_t$ is analytic at (0,0). 
Imitate the former proof noting that, so defined, $\nu_t$ is analytic at (0,0).\\
{\underline {${\cal C}^\infty$ case :}} Let $\tilde{T_t}=j_0^\infty(T_t)$ and
$\tilde{\varepsilon_t}\in\R_t\big[[x,y]\big]$ be such that 
$W.\tilde{\varepsilon_t}=\tilde{T_t}$. 
Borel's theorem ensures 
the existence of $\varepsilon_t\in{\cal C}^\infty(\R^2)$ such that 
$j_0^\infty(\varepsilon_t)=\tilde{\varepsilon_t}$. 
Thus $W.\varepsilon_t=T_t+m_t$ 
where $m_t$ is flat at (0,0). 
Let $n_t$ be flat at (0,0) and such that $W.n_t=-m_t$
($n_t$ exists by theorem \ref{roussarie}); $\nu_t=\varepsilon_t+n_t$
suits.
$\qed$\\

\noindent Consequently, to prove {\bf 1-}, we put $\beta_t=exp\,\nu_t$, where
$\nu_t$ satisfies $W.\nu_t=\mu_t$.\\

{\bf 2-} Note first that if $\degp$ is a quasidegree (for polynomials), then 
there exists $(i,j)$ in $\N^2$ such that $\degp=i\omega_1+j\omega_2$
if not, $\degp=i\omega_1-\omega_2$ (or $-\omega_1+i\omega_2$) with
$i\in\N$. 
The following lemma will prove the second claim.
\begin{lem}  
\begin{description}
\item[i)] Let $k$ and $l$ be in $\N$ and $T_t\in{\cal F}_t(\K^2)$ with 
${\rm {ord}}\big((j_0^\infty(T_t)\big)>k\omega_1+l\omega_2$. 
Then there exists 
$\gamma_t\in{\cal F}_t(\K^2)$
satisfying $W.\gamma_t-(k\omega_1+l\omega_2)\gamma_t=T_t$.
\item[ii)] Let $k\in\N$ and $T_t\in{\cal F}_t(\K^2)$; then there exists 
$\gamma_t\in{\cal F}_t(\K^2)$ satisfying \\
$W.\gamma_t-(k\omega_1-\omega_2)\gamma_t=T_t$.
\end{description}
\label{E2}
\end{lem}   
{\sf Proof }: 
{\bf i)} We use an induction :\\
For $k=l=0$: see lemma \ref{E1}.\\
Now, assume that i) is true for $(k,l)\in\N^2$. 
We are going to show that it is
true for $k+1$ and $l$ (for $k$ and $l+1$ the proof is the same). \\
Let $T_t\in{\cal F}_t(\K^2)$ with 
${\mbox {ord}}\big(j_0^\infty(T_t)\big)>(k+1)\omega_1+l\omega_2$ and 
$\delta_t\in{\cal F}_t(\K^2)$ verifying 
$W.\delta_t-(k\omega_1+l\omega_2)\delta_t=\frac{\partial T_t}{\partial x}$.
Then we define $\gamma_t$ by $\gamma_t(x,y)=\int_0^x\delta_t(u,y)\, du$ for 
$(x,y)$ in a neighbourhood of (0,0). 
An easy computation shows that
$W.\gamma_t-((k+1)\omega_1+l\omega_2)\gamma_t=T_t$.\\
{\bf ii)} We use again an induction :\\
For $k=0$: we know that there exists $\delta_t\in{\cal F}_t(\K^2)$ such that 
$W.\delta_t=\int_0^yT_t(x,u)\, du$. 
If we put 
$\gamma_t=\frac{\partial \delta_t}{\partial y}$ then we get
 $W.\gamma_t+\omega_2\gamma_t=T_t$.\\
The end of the proof can be achieved as in i).
$\qed$\\

\noindent {\bf A list of normal forms :} 
We recall that a germ of Poisson structures on $\K^n$ is determined by the choice of
a germ of functions.\\
We consider a germ $\Pi=f\biv$, where $f$ vanishes and is singular at $(0,0)$.
We suppose, in addition, that the germ $f$ is of {\bf finite codimension}.
It means that the vector space $Q_f={\cal F}(\K^2)/I_f$ ($I_f$ is the ideal spanned by 
$\frac{\partial f}{\partial x}$ and $\frac{\partial f}{\partial y}$)
is of finite dimension.\\
Why do we suppose that $f$ is of finite codimension? In fact, one can see $I_f$ as the
tangent space of the orbit of $f$ (with respect to the R-equivalence).
 Thus, the 
finite-codimensional germs are those whose orbit is ``big'' enough.

\begin{rem}
{\rm It is important to note the following fact:\\
According to Tougeron's theorem (see for instance [AGV]), if $f$ is of finite 
codimension, then $f$ is R-equivalent to its $k$-jets for $k$ sufficiently large.
The set $f^{-1}(\{0\})$ is then, from the topological 
point of view, the same as the set of zeroes of a polynomial.
Therefore, if $g$ is a germ at 0 of functions which satisfies $fg=0$, then $g=0$.}
\label{codim}
\end{rem}
\endrem

Moreover, we suppose that the germ $f$ is {\bf simple}. 
It means
that a sufficiently small neighbourhood (with respect to Whitney's topology; see [AGV]) 
of $f$ intersects only a finite number of R-orbits. 
Simple germs are those who present 
a certain kind of stability under deformation.\\
Note that simple germs are necessarily of finite codimension (see for
instance [AGV]).\\
We have a classification of such germs in the following theorem.
\begin{theo} {\rm {\bf [AGV]}}
Simple germs at (0,0) of functions are given, up to R-equivalence, in the following list: 
\begin{center}
\begin{tabular}{|c|c|c|c|c|}
      \hline
      $A_k \quad k\geq 1$ & $D_k \quad k\geq 4$ & $E_6$ & $E_7$ & $E_8$\\
      \hline
      $x^2\pm y^{k+1}$ & $x^2y \pm y^{k-1}$ & $x^3\pm y^4$ & $x^3 + xy^3$ & $x^3 + y^5$\\
      \hline
      \end{tabular}\\
\end{center}
If $\K=\C$ (or if $k$ is even in the real $A_k$ case), then the symbol $\pm$ disappears.
\label{germes simples}
\end{theo}
It is important to note that these models are \qh polynomials.\\
Now, applying theorem \ref{fn poisson} to these models, we can state the following theorem.
\begin{theo} {\rm {\bf [A]}}
Let $f$ be a simple germ at (0,0) of finite codimension. 
Suppose that $f$ has 
at (0,0) a critical point with critical value 0. 
Then, the germ $\Pi=f\biv$ 
is equivalent, up to a multiplicative constant, to a germ of type 
$g\biv$, where $g$ is in the following list:
\begin{eqnarray*}
A_{2p} &:& x^2+y^{2p+1}\quad p\geq 1\\
A_{2p-1}^{\pm} &:& (x^2\pm y^{2p})(1+\lambda y^{p-1})\quad p\geq 1\\
D_{2p}^{\pm} &:& (x^2\pm y^{2p})(1+\lambda y^{p-1})\quad p\geq 2\\
D_{2p+1} &:& (x^2y+y^{2p})(1+\lambda x)\quad p\geq 2\\
E_6 &:& x^3+y^4\\
E_7 &:& (x^3+xy^3)(1+\lambda y^2)\\
E_8 &:& x^3+y^5
\end{eqnarray*}
If $\K=\C$, the symbol $\pm$ disappears.
\label{modeles}
\end{theo}

\section{Poisson cohomology}
In this section we compute the Poisson cohomology of some Poisson structures.
In fact, we work {\bf locally} and we study the {\it ``germified''}
Poisson cohomology. 
This means 
that we work with germs of Poisson structures, functions, vector fields, 2-vectors...\\
We recall that $\fct$ ($\chp$, $\pois$) indicates the space of germs at 0 of
functions (vector fields, 2-vectors).\\
Let $\Pi$ be a germ of Poisson structure on $\K^2$. 
We have then the complex
$$
0\stackrel{\delta_0}{\longrightarrow}\fct\stackrel{\delta_1}{\longrightarrow}\chp
\stackrel{\delta_2}{\longrightarrow}\pois\stackrel{\delta_3}{\longrightarrow}0
$$
where $\delta_0=0$, $\delta_3=0$, $\delta_1(g)=[g,\Pi]$ and 
$\delta_2(X)=[X,\Pi]$ ([.,.] indicates the Schouten bracket).\\
We will denote by $Z^i(\Pi)=\mbox{Ker }\delta_{i+1}$, 
$B^i(\Pi)=\mbox{Im }\delta_i$ and 
$H^i(\Pi)=Z^i(\Pi)/B^i(\Pi)$.\\
If we assume that $\Pi=F\biv$ where $F\in\fct$ then, for $g\in\fct\!$
, we have 
$$
\delta_1(g)=F\frac{\partial g}{\partial y}\vx-
F\frac{\partial g}{\partial x}\vy\,.
$$ 
We will denote by $X_g$ this vector field 
(it is the Hamiltonian of $g$ with respect to $\Pi$) and $H_g$ the vector field
$\frac{\partial g}{\partial y}\vx-\frac{\partial g}{\partial x}\vy$.\\
On the other hand, for $X\in\chp$, we have 
$$
\delta_2(X)=\big(X.F-(div\!X)F\big)\biv\,.
$$ 
We will denote $H^2(F)$ the space 
$\fct/\{X.F-(div\!X)F\| X\in\chp\}$. 
This space is clearly isomorphic to $H^2(\Pi)$.\\

Actually, we will compute the cohomology of Poisson structures of a particular type.\\
Let $(\omega_1,\omega_2)\in\N^*\times\N^*$. 
Here and throughout, the quasihomogeneity will be understood as being in the sense of
$(\omega_1,\omega_2)$ ($W$ will again indicate the vector field 
$\omega_1x\vx+\omega_2y\vy$).\\
We take a \qh polynomial $f$ of degree d and we assume that $f$ is a germ at 0 of 
finite codimension $c$ (recall that it means that the vector space $Q_f=\fct/I_f$, where
$I_f$ is the ideal spanned by $\frac{\partial f}{\partial x}$ and  
$\frac{\partial f}{\partial y}$, is of dimension $c$). 
We also give us a \qh polynomial
$h$ of degree $\degp$ (if $\degp>0$).\\
Now we consider two germs of Poisson structures
$$
\Pi_0=f\biv\quad {\mbox { and }}\quad \Pi=f(1+h)\biv
$$
and we are going to compute the cohomology of these Poisson structures.\\
In the former section, we saw that the ``most interesting'' Poisson structures are of 
this type.\\
Note that {\bf {in the sequel,we do not suppose that the germ $f$ is
simple}}.\\

In our computation, it is very important to assume that $f$ is of finite
codimension (see the role played by the second claim of lemma \ref{l1}, 
and the remark \ref{Rk}).\\ 
It is easy to see that, since $f$ is of finite codimesion, the spaces $H^0(\Pi_0)$ and
$H^0(\Pi)$ are isomorphic to $\K$ (see remark \ref{codim}).
\begin{rem}
{\rm It is important to note that, since $f$ is quasihomogeneous, the computation of 
$H^\bullet(\Pi_0)$ 
can be done "degree by degree". 
For instance, if $\sum_i X^{(i)}$ is the $\infty$-jet of $X$ and
if $X$ is in $Z^1(\Pi_0)$ (resp. $B^1(\Pi_0)$) then $X^{(i)}$ is also in $Z^1(\Pi_0)$ 
(resp. $B^1(\Pi_0)$) for each $i$. 
Moreover, if $X$ is polynomial, then $X$ is in $Z^1(\Pi_0)$
(resp. $B^1(\Pi_0)$) if and only if each of its \qh components is in $Z^1(\Pi_0)$ (resp.
$B^1(\Pi_0)$). 
We have the same properties for $B^2(\Pi_0)$.\\
The computation of the cohomology of $\Pi$ does not present this property.}
\label{dad}
\end{rem}
\endrem
The following result will be useful in the sequel.
\begin{lem}
Let $X$ be in $\chp$.\\
1- If $div\,X=0$, then there exists $g\in\fct$ such that $X=H_g$.\\
2- If $X.f=0$, then $X=\alpha H_f$ with $\alpha \in \fct$.
\label{l1}
\end{lem}
{\sf Proof }: We can write $X=A \vx +B\vy$. 
We consider the 1-form $\omega=-Bdx+Ady$.\\
1- If $div\,X=0$  then $d\omega=0$, which implies that $\omega=dg$ with $g\in\fct$, and so
$X=H_g$.\\
2- If $X.f=0$ then $df\wedge \omega=0$. 
Since $f$
has finite codimension, de Rham's division theorem (see [dR] or [M]) enables us to 
conclude.
$\qed$
\subsection{Computation of $H^1$}
{\underline {{\bf Computation of $H^1(\Pi_0)$ :}}}
\begin{lem}
Let $X\in Z^1(\Pi_0)$. 
Then there exists $\alpha\in\fct$ such that\\
$X=\alpha H_f+\frac{div\!X}{{\rm d}}W$.
\label{l2}
\end{lem}
{\sf Proof }: Direct application of lemma \ref{l1}.\\

The main idea in the computation of the space $H^1(\Pi_0)$, is to show that every 1-cocycle
whose $\infty$-jet has a sufficiently large order is a cobord. 
\begin{lem}
Let $X\in\chp$ be such that ${\rm {ord}}\big(j_0^\infty(X)\big)>\degp$.\\
If $X\in Z^1(\Pi_0)$ then $X\in B^1(\Pi_0)$.
\label{l3}
\end{lem}
{\sf Proof }: {\underline {first case : $div\,X=0$.}} 
We then show that $f$
divides X.\\
Since $X.f-(div\,X)f=0$, we have $X.f=0$ and then $X=\gamma H_f$ (lemma \ref{l1}) with 
$\gamma\in\fct$. 
Note that, since ${\mbox {ord}}\big(j_0^\infty(X)\big)>\degp$ and $H_f$ 
is \qh
of degree $\degp$, we have ${\mbox {ord}}\big(j_0^\infty(\gamma)\big)>0$.\\
We prove that $f$ divides $\gamma$.
Let $\mu\in\fct$ be such that $W.\mu=\gamma$ ($\mu$ exists according to lemma 
\ref{E1}). 
Note that, if we write the $\infty$-jet in the relation $W.\mu=\gamma$, we see
that the order of the $\infty$-jet of $\mu$ is also strictly positive 
(because $W$ is \qh of degree 0). 
Therefore, the order of the $\infty$-jet of $H_f.\mu$
is strictly larger than $\degp$, because $H_f$ is \qh of degree $\degp$.\\ 
We have 
$$
H_f.(W.\mu)=W.(H_f.\mu)+[H_f,W].\mu=W.(H_f.\mu)+\big(-(\degp)H_f\big).\mu
$$
(because $H_f$ is of degree $\degp$).\\
Since $H_f.(W.\mu)=H_f.\gamma=div\,X=0$,
we have 
$$
W.(H_f.\mu)=(\degp)H_f.\mu
$$ 
and so $H_f.\mu$ is either 0 or \qh of degree $\degp$. \\
Thus, as 
${\mbox {ord}}\big(j_0^\infty(H_f.\mu)\big)>\degp$, we have $H_f.\mu=0$.\\
Now, $H_\mu.f=-H_f.\mu=0$ so there exists $\nu\in\fct$ such that
$\frac{\partial\mu}{\partial x}=\nu\frac{\partial f}{\partial x}$ and
$\frac{\partial\mu}{\partial y}=\nu\frac{\partial f}{\partial y}$
(lemma \ref{l1}). 
Therefore, $W.\mu=\nu\, W.f$, that is
$\gamma=\nu ({\rm d}\times f)$.\\
We deduce that $X=fZ$ with $Z\in\chp$.\\ 
Finally, since $X\in Z^1(\Pi_0)$, 
$div\,Z=0$ and then $Z=H_g$ for some $g\in\fct$ (lemma \ref{l1}). 
Hence $X=fH_g=X_g$.\\

\noindent {\underline {Second case : $div\,X\neq0$.}}
If we find $\beta\in\fct$ such that $div\,X=div\,X_\beta$, then the 1-cocycle $X-X_\beta$
satisfies $div(X-X_\beta)=0$, which implies (see the first case) that 
$X=X_\beta+X_\varepsilon$ where $\varepsilon\in\fct$. 
Since 
$div\,X_\beta=H_\beta.f=-H_f.\beta$, we are looking for $\beta$ such that 
$H_f.\beta=-div\,X$.\\
We have $X=\alpha H_f+\frac{div\,X}{{\rm d}}W$ with $\alpha\in\fct$ 
(lemma \ref{l2}) so that $div\,X$ satisfies the equation 
$$
W.(div\,X)-(\degp)div\,X=-{\rm d}\times H_f.\alpha\,.
$$ 
Note that, if we write the $\infty$-jets at 0 in the relation 
$X=\alpha H_f+\frac{div\,X}{{\rm d}}W$, we see that the order of the $\infty$-jet of 
$\alpha$ is strictly positive so, according to lemma \ref{E1} we can take $\beta\in\fct$ 
such that $W.\beta={\rm d}\times \alpha$ (the order of
the $\infty$-jet of $\beta$ is also strictly positive). 
Since
$$
W.(H_f.\beta)=H_f.(W.\beta)+[W,H_f].\beta\\
={\rm d} (H_f.\alpha)+\big((\degp)H_f\big).\beta\,,
$$
we have $\quad W.(div\,X+H_f.\beta)=(\degp)(div\,X+H_f.\beta)$.\\
We deduce that $div\,X+H_f.\beta$ is either 0 or \qh of degree $\degp$. 
Therefore, $div\,X=-H_f.\beta$ 
(because ${\mbox {ord}}\big(j_0^\infty(div\,X+H_f.\beta)\big)>\degp$).
$\qed$\\

\noindent For $X\in Z^1(\Pi_0)$, we denote by ${[X]}_{ _{\Pi_0}}$ its class 
modulo $B^1(\Pi_0)$. 
We also denote by $\{e_1,...,e_r\}$ a basis of the vector
space of \qh polynomials of degree $\degp$ (in fact, in order to
obtain a vector space, we have to add 0 to the set of \qh 
polynomials of degree $\degp$).
\begin{theo}
The family $\big\{ [H_f]_{ _{\Pi_0}},[e_1W]_{ _{\Pi_0}},...,[e_rW]_{ _{\Pi_0}} \big\}$ 
is a basis of $H^1(\Pi_0)$. 
In particular, $H^1(\Pi_0)$ is a finite-dimensional vector space
of dimension $r+1$.
\label{t1}
\end{theo}
{\sf Proof }: First we prove that $H^1(\Pi_0)$ is spanned by this family.
Lemma \ref{l3} says that every $X$ in $Z^1(\Pi_0)$ is cohomologous to a polynomial
vector field of maximum degree $\degp$. 
Indeed, if $X\in Z^1(\Pi_0)$ then $j_0^\degp(X)$ is also
in $Z^1(\Pi_0)$ ($j_0^\degp(X)$ indicates the jet of degree $\degp$ of
$X$ at 0). 
Thus, 
$X-j_0^\degp(X)$ is in $Z^1(\Pi_0)$ and the order of its $\infty$-jet at 0 is strictly higher
than $\degp$.\\ 
Therefore, using remark \ref{dad}, we can assume that $X$ is \qh of degree 
$\degp$ or lower.\\
- If $X\in Z^1(\Pi_0)$ is \qh with $\deg X<\degp$ then $X=0$. 
Indeed, according to
lemma \ref{l2}, we have $X=\frac{div\,X}{{\rm d}}W$, and so 
$$
div\,X=\frac{div\,X}{{\rm d}}div\,W+W.\big(\frac{div\,X}{d}\big)\,,
$$ 
which implies that
$(\degp-\deg X)div\,X=0$.\\
- Let $X\in Z^1(\Pi_0)$ be quasihomogeneous of degree $\degp$. 
We have (lemma \ref{l2})
$X=\alpha H_f+\frac{div\,X}{{\rm d}}W$ where $\alpha\in\K$ and $div\,X$ is a \qh
polynomial of degree $\degp$.\\ 
Therefore the family generates $H^1(\Pi_0)$.\\
Now, we prove that this family is free.
Suppose that  $\sum_i\lambda_ie_iW+\alpha H_f\in B^1(\Pi)$
where $\alpha,\lambda_1,\,.\,.\,.\,, \lambda_r$
are scalars. 
Then $\sum_i\lambda_ie_iW+\alpha H_f=0$. 
Indeed, if $g$ is a \qh
polynomial, then $\deg X_g=\deg H_f+\deg g=\degp + \deg g$, which is strictly larger
than $\degp$ as soon as $g\neq 0$.\\
Consequently, $div\big(\sum_i\lambda_ie_iW+\alpha H_f\big)=0$ i.e. 
$\sum_i\lambda_ie_i=0$. 
We deduce that $\lambda_1,\hdots,\lambda_r$
are 0, and so $\alpha=0$.
$\qed$\\

\noindent{\underline {{\bf Computation of $H^1(\Pi)$ :}}}\\
If $X\in Z^1(\Pi)$ we denote by $[X]_{ _{\Pi}}$ its class modulo $B^1(\Pi)$.
\begin{theo}
$\big\{ [(1+h)H_f]_{ _{\Pi}},[(1+h)e_1W]_{ _{\Pi}},...,
[(1+h)e_rW]_{ _{\Pi}}\big\}$ is a 
basis of $H^1(\Pi)$. 
In particular, $H^1(\Pi)\simeq H^1(\Pi_0)$.
\end{theo}
{\sf Proof }: It is sufficient to notice that 
$X\in Z^1(\Pi)$ (resp. $X\in B^1(\Pi)$) if and only if $\frac{X}{1+h}\in Z^1(\Pi_0)$
(resp. $B^1(\Pi_0)$).
$\qed$

\subsection{Computation of $H^2$}
{\underline {{\bf Computation of $H^2(\Pi_0)$ :}}}
\begin{lem}
Let $g$ be a germ at 0 of functions on $\K^2$.\\
1. If the $\infty$-jet at $0$ of $g$ does not contain a component of degree $2{\rm d}-w_1-w_2$ 
then 
$$
g\in B^2(f)\Leftrightarrow g\in I_f\,.
$$
2. If $g$ is \qh of degree $2{\rm d}-w_1-w_2$ then  
$$
g\in B^2(f)\Rightarrow g\in I_f\,.
$$ 
\label{l4}
\end{lem}
{\sf Proof }:
If $g=X.f-(div\,X)f\in B^2(f)$ where $X\in\chp$ then , if we put 
$$
Y=X-\frac{div\,X}{{\rm d}}W\,,
$$ 
we have $g=Y.f$\\
This proves the second claim and the first part of the first one.\\
Now, we prove the converse of the first claim: we assume that $g\in I_f$ (and the $\infty$-jet 
of $g$ does not contain a component of degree $2{\rm d}-w_1-w_2$) and we are going to
show that $g\in B^2(f)$.\\
{\it Formal case }: Let $\D g=\sum_{i\geq 0}g^{(i)}\in\K[[x,y]]$ and
$\D X=\sum_{i\geq {\rm d}-\max (\omega_1,\omega_2)}X^{(i-{\rm d})}$ 
(with $g^{(i)}$ of degree $i$ and $X^{(i-{\rm d})}$ of degree $i-{\rm d}$) such that 
$g=X.f$. 
Note that $X^{(\degp)}=0$.\\ 
If we put 
$$
Y=X+\sum_{i\neq 2{\rm d}-\omega_1-\omega_2}\frac{div\, X^{(i-{\rm d})}}
{2d-\omega_1-\omega_2-i}\,W\,,
$$ 
a direct computation gives $Y.f-(div\,Y)f=X.f=g$.\\
{\it Analytical case }: If $X$ is analytic at (0,0), then $div\, X$ is analytic too and
since \\ $lim_{i\rightarrow +\infty} \frac{1}{2{\rm d}-\omega_1-\omega_2-i}=0$ the 
vector field defined above is also analytic in (0,0).\\
{\it ${\cal C}^\infty$ case }: Let us denote by ${\tilde {g}}=j_0^\infty (g)$ and
${\tilde {X}}=j_0^\infty (X)$.\\ 
If we write the $\infty$-jets in the relation $g=X.f$, we get ${\tilde {g}}={\tilde {X}}.f$. 
Thus, there exists a formal vector field ${\tilde {Y}}$ such that
$$
{\tilde {g}}={\tilde {Y}}.f-(div\, {\tilde {Y}})f\,.
$$
Let $Y$ be a ${\cal C}^\infty$ vector field such that ${\tilde {Y}}=j_0^\infty(Y)$. 
Since $Y.f-(div\, Y)f$ and $g$ have the same $\infty$-jet, this vector field satisfies 
$$
Y.f-(div\, Y)f=g+\varepsilon
$$ 
where $\varepsilon$ 
is flat at 0.\\
Now, since $Y.f-(div\, Y)f\in I_f$ (see the beginning of the proof), $\varepsilon$ is in 
$I_f$ so that $\varepsilon=P.f$ where $P$ is a flat vector field. 
According to lemma \ref{E2}, 
there exists $\alpha\in\fct$ such that $W.\alpha -(\degp)\alpha =-div\, P$.
Consequently, setting $Z=P+\alpha W$, we have $Z.f-(div\,
Z)f=\varepsilon$.
$\qed$\\

\begin{rem}
{\rm 1- This lemma is true even if $f$ is not of finite codimension.\\
2- This lemma gives $B^2(f)\subset I_f$. 
Thus, there is a surjection from $H^2(f)$ onto
$Q_f$. 
Therefore, if $f$ is not of finite codimension, then $H^2(\Pi_0)$ is an 
infinite-dimensional vector space.\\
3- Finally, according to this lemma, if $\xi\in I_f$, then there exists a \qh polynomial 
${\overline {\xi}}$ of degree $2{\rm d}-\omega_1-\omega_2$ such that 
$\xi+{\overline {\xi}}\in B^2(f)$.}
\label{Rk}
\end{rem}
\endrem

\noindent If $g\in\fct$, $[g]_{ _{\Pi_0}}$ indicates its class modulo $B^2(f)$.
Recall that $\{e_1,\hdots,e_r\}$ is a basis of the space of \qh polynomials of degree $\degp$.
Finally, we denote by $\{u_1,...,u_c\}$ a monomial
basis of $Q_f=\fct/I_f$ (for the existence of such a basis, see [AGV]).
\begin{theo}
The family
$\big\{ [e_1f]_{ _{\Pi_0}},...,[e_rf]_{ _{\Pi_0}},[u_1]_{ _{\Pi_0}},...,
[u_c]_{ _{\Pi_0}}\big\}$ is a basis of $H^2(f)$.\\
In particular, $H^2(\Pi_0)$ is a finite-dimensional vector space of dimension $r+c$.
\label{t3}
\end{theo}
{\sf Proof }:- This family generates $H^2(f)$:\\
Let $g\in\fct$. 
We can write $g=\sum_{i=1}^c \lambda_i u_i+\xi$ where $\lambda_i\in\K$ and
$\xi\in I_f$. 
According to lemma \ref{l4} (and remark \ref{Rk}), we can write
$$
g=\sum_{i=1}^c \lambda_i u_i + {\overline {g}}\quad {\mbox
{mod}}\,B^2(f)
$$
where ${\overline g}$ is a \qh polynomial of degree $2{\rm d}-\omega_1-\omega_2$.\\
We can again write ${\overline {g}}=\sum_{i=1}^c {\overline {\lambda_i}} u_i\, {\mbox { mod }}
\, I_f$ where ${\overline {\lambda_i}}\in\K$.
Now, we know (see [AGV] p.200) that 
$\max \{\deg u_1,...,\deg u_c\}=2{\rm d}-2\omega_1-2\omega_2$
which is strictly lower than $\deg {\overline {g}}$. 
So, ${\overline {g}}$ is in $I_f$, i.e.
${\overline {g}}=X.f$ with $X$ \qh of degree $\degp$.\\
Therefore, ${\overline g}=(div\, X)f+\big( X.f-(div\, X)f\big)$ with $div\, X$ \qh of 
degree $\degp$.\\
- This family is free: Let $g_1=\sum_{i=1}^r \lambda_i e_i$ and 
$g_2=\sum_{j=1}^c \mu_j u_j$ with $\lambda_i$ and $\mu_j$ in $\K$ for any $i$
and $j$. 
We assume that $g_1f+g_2\in B^2(f)$. 
Since
$\max \{\deg u_1,...,\deg u_c\}<\deg (g_1f)$, $g_1f$ and $g_2$ are both in $B^2(f)$ 
(see remark \ref{dad}). \\
On the one hand, we then have $g_2\in I_f$ which is possible only if 
$\mu_1=...=\mu_c=0$. \\
On the other hand, since $g_1f=X.f-(div\, X)f$ for some quasihomogeneous vector field $X$
of degree $\degp$, if $Y$ denotes the vector field $\frac{g_1+div\, X}{{\rm d}}W$
then $(X-Y).f=0$. 
Therefore (lemma \ref{l1}), $X=Y+\alpha H_f$ with $\alpha\in\K$.
Hence $X.f-(div\, X)f=0$, which implies $\lambda_1=...=\lambda_r=0$.
$\qed$\\

\noindent {\underline {{\bf Computation of $H^2(\Pi)$ :}}}
\begin{lem}
Let $g\in\fct$. 
If the order of the $\infty$-jet of $g$ is larger than or equal to
$2{\rm d}-\omega_1-\omega_2$, then there exists a \qh polynomial $\varepsilon$
of degree $\degp$ such that
$g=\varepsilon f \quad {\rm {mod}}\, B^2(f+fh)$.
\label{l5}
\end{lem} 
{\sf Proof }: 
According to theorem \ref{t3}, we can write
$$
\frac{g}{1+h}= \sum_{i=1}^c \lambda_iu_i+ \varepsilon f+X.f-(div\, X)f
$$
where the $\lambda_i$ are in $\K$ and $\varepsilon$ is \qh of degree $\degp$.\\
But since ${\mbox {ord}}\Big( j_0^\infty\big(\frac{g}{1+h}\big)\Big)\geq 2{\rm d}-
\omega_1-\omega_2>\max \{\deg u_1,\hdots,\deg u_c\}$, if we write the $\infty$-jets in 
the former relation, we see that the $\lambda_i$ are zero and that the 
order of the $\infty$-jet of $X$ is higher than $\degp$.\\
Now, 
$$
X.(f+fh)-(div\, X)(f+fh)=g+f(X.h)-\varepsilon f(1+h)\,.
$$
We put $\lambda=(\degp)\big( 1+\frac{h}{1+h}\big)$. 
If $X.h-\varepsilon h=0$, then the lemma
is shown. 
Now, if we suppose that $X.h-\varepsilon h$ is not 0, the order of its $\infty$-jet
is strictly larger than ${\rm d}-\omega_1-\omega_2$, and so we can take $\alpha\in\fct$ 
such that $W.\alpha -\lambda\alpha=\frac{X.h-\varepsilon h}{1+h}$ (see the  "Resolubility of 
equation $(E)$" in the former section). \\
If we put $Z=X+\alpha W$, we have $Z.(f+fh)-(div\,
Z)(f+fh)+f\varepsilon=g$.
$\qed$
\begin{theo}
The family $\big\{ [e_1f]_{ _{\Pi}},...,[e_rf]_{ _{\Pi}},[u_1]_{ _{\Pi}},...,
[u_c]_{ _{\Pi}}\big\}$ is a basis of $H^2(f+fh)$. 
In particular, 
$H^2(f+fh)\simeq H^2(f)$ (the space $H^2(\Pi)$ is then of dimension $r+c$).
\label{t4}
\end{theo}
{\sf Proof }:\\
- This family generates $H^2(f+fh)$.\\ 
Given $g\in\fct$, we have 
$g=\sum_{i=1}^c \lambda_{i,0}u_i+P_0f+X_0.f-(div\, X_0)f$ (theorem \ref{t3}) where 
$P_0$ is a \qh polynomial of degree $\degp$, $\lambda_{i,0}\in\K$ for any $i$ and
${\mbox {ord}}\big( j_0^\infty (X_0)\big)\geq -{\mbox {max}}(\omega_1,\omega_2)$
(remember that $-{\mbox {max}}(\omega_1,\omega_2)$ is the smallest degree for \qh
vector fields).
Since
$$
X_0.(f+fh)-(div\, X_0)(f+fh)=g-\sum_{i=1}^c \lambda_{i,0} u_i-P_0f+
X_0.(fh)-(div\, X_0)(fh)\; ,
$$
we can write
$$
g=\sum_{i=1}^c \lambda_{i,0}u_i +P_0f-\big(X_0.(fh)-
(div\, X_0)(fh)\big)\quad {\mbox {mod}}\,B^2(f+fh).
$$
Now, we have $X_0.(fh)-(div\, X_0)fh=-\sum\lambda_{i,1}u_i-P_1f+X_1.f-(div\, X_1)f$
where $\lambda_{i,1}$ is in $\K$ for any $i$,
$P_1$ is a \qh polynomial of degree $\degp$ 
and ${\mbox {ord}}\big( j_0^\infty (X_1)\big)\geq \degp-{\mbox {max}}
(\omega_1,\omega_2)$. 
So, in the same way, 
$$
X_0.(fh)-(div\, X_0)fh=-\sum_{i=1}^c \lambda_{i,1}u_i-P_1f-\big(X_1.(fh)-
(div\, X_1)(fh)\big) \;{\mbox {mod}}\, B^2(f+fh).
$$
Hence 
$$
g=\sum_{i=1}^c (\lambda_{i,0}+\lambda_{i,1})u_i+
(P_0+P_1)f-\big(X_1.(fh)-(div\, X_1)(fh) \big) \quad {\mbox {mod}}\, 
B^2(f+fh).
$$
In this way, we get
$$
g=\sum_{i=1}^c (\lambda_{i,0}+...+\lambda_{i,k})u_i+
(P_0+...+P_k)f-\big(X_k.(fh)-(div\, X_k)(fh)\big)\;{\mbox {mod}}\,
B^2(f+fh)
$$
where $k$ is the smallest integer such that $k(\degp)-
{\mbox {max}}(\omega_1,\omega_2)\geq 0$,
$P_j$ is a \qh polynomial of degree $\degp$ for any $j$, 
$\lambda_{i,j}\in\K$ for any $i$ and $j$ and
${\mbox {ord}}\big( j_0^\infty (X_k)\big)\geq k(\degp)-{\mbox {max}}
(\omega_1,\omega_2)$.\\
But since 
\begin{eqnarray*}
\ {\mbox {ord}}\big( j_0^\infty (X_k.(fh)-(div\, X_k)fh)\big) & \geq & 
2{\rm d}-\omega_1-\omega_2+k(\degp)- {\mbox {max}}(\omega_1,\omega_2)\\
\ & \geq & 2{\rm d}-\omega_1-\omega_2\: ,
\end{eqnarray*}
 lemma \ref{l5} gives $X_k-(div\, X_k)f=Qf \quad {\mbox {mod}}\, 
B^2(f+fh)$ for some \qh polynomial $Q$ of degree $\degp$.\\
- This family is free. 
Let $\lambda_1,...\lambda_r$ be scalars and $P$ a \qh
polynomial of degree $\degp$. 
Suppose that 
$$
\sum_{i=1}^c \lambda_iu_i+Pf=X.(f+fh)-(div\, X)(f+fh)\quad (\ast)
$$ 
with $X\in\chp$. 
Since
$fh$ is \qh of degree $2{\rm d}-\omega_1-\omega_2$ and the order of the $\infty$-jet 
of $X$ is larger than $-{\mbox {max}}(\omega_1,\omega_2)$, we have
\begin{eqnarray*}
\ {\mbox {ord}}\big( j_0^\infty (X.(fh)-(div\, X)(fh))\big) & \geq &
 2{\rm d}-\omega_1-\omega_2- {\mbox {max}}(\omega_1,\omega_2)\\
\ & > & 2(\degp)={\mbox {max}}\{{\mbox {deg}}u_1,...,{\mbox {deg}}u_c\}\: .
\end{eqnarray*}
Therefore, if we write the $\infty$-jets in the relation $(\ast)$, we have 
$\sum_{i=1}^c\lambda_iu_i\in B^2(f)$ and so 
$\lambda_1=\,,\,.\,.\,.\,,=\lambda_c=0$ (theorem \ref{t3}).\\
We obtain 
$$
Pf=X.(f+fh)-(div\, X)(f+fh)\quad (\ast\ast)\,,
$$ 
where $Pf$ is a 
\qh polynomial of degree $2{\rm d}-\omega_1-\omega_2$.\\
Now, we can write $j_0^\infty(X)=\sum_{i\geq \delta} X^{(i)}$ ($X^{(i)}$ is \qh of
degree $i$). 
If $\delta <\degp$ then $X^{(\delta)}\in Z^1(\Pi_0)$ and so, 
$X^{(\delta)}=0$ (cf proof of theorem \ref{t1}).\\
In the same way, we can prove that $X^{(i)}=0$ for any $i<\degp$.\\
Consequently, if we write the $\infty$-jets in the relation $(\ast\ast)$, we get
$$
Pf=X^{(\degp)}.f-div\,X^{(\degp)}f
$$ 
that is, $Pf\in B^2(f)$,
which is possible only if $P=0$ (cf proof of theorem \ref{t3}).
$\qed$

\subsection{Examples}
We are going to make explicit the cohomology of some Poisson structures given in theorem 
\ref{modeles}.\\

\noindent {\it The regular case} : We suppose that $\Pi=x\biv$.\\
In this case, $\omega_1=\omega_2={\rm d}=1$, and so $\degp<0$.
Moreover, $Q_x=\{0\}$.\\
Therefore, $H^1(\Pi)\simeq \K \vy$ and $H^2(\Pi)=\{0\}$.\\

\noindent {\it Morse's singularity \; $(A_1)$}: We suppose that $\Pi=(x^2+y^2)\biv$.\\
Here, we have $\omega_1=\omega_2=1$ and ${\rm d}=2$. 
The only monomials of degree 
${\rm d}-\omega_1-\omega_2$ are the scalars. 
Moreover, $Q_{x^2+y^2}\simeq\K.1$. 
Then, we have 
\begin{eqnarray*}
H^1(\Pi) &\simeq& \K.(y\vx-x\vy)\oplus\K.(x\vx+y\vy)\\
{\mbox { and }} H^2(\Pi) &\simeq& \K.\biv\oplus \K.f\biv.
\end{eqnarray*}

\noindent {\it The singularity $D_{2p+1}\; (p\geq 2)$}: We suppose that
$\Pi=(x^2y+y^{2p})(1+x)\biv$.\\
In this case, we can see that $\omega_1=2p-1$, $\omega_2=2$ and ${\rm d}=4p$. 
The only monomials of degree ${\rm d}-\omega_1-\omega_2$ are of type $\lambda x$ 
(with $\lambda\in\K$).
Moreover, the family $\{ 1,x,y,y^2,\cdots,y^{2p} \}$ is a monomial basis of 
$Q_{x^2y+y^{2p}}$.\\
Therefore,\\ 
the family $\Big\{ \big[ (1+x)\big((x^2+2py^{2p-1})\vx-2xy\vy\big) \big],
\big[ (1+x)x\,W \big] \Big\}$ is a basis of $H^1(\Pi)$\\ 
and the family \\ 
$\Big\{ \big[ x(x^2y+y^{2p})\biv\big],\big[\biv\big],\big[x\biv\big],\big[y\biv\big],
\big[y^2\biv\big],\,.\,.\,.\,,$\\
$\,.\,.\,.\,,\big[y^{2p}\biv\big]\Big\}$ is a basis of $H^2(\Pi)$.\\
In particular, $H^1(\Pi)$ is of dimension $2$ and $H^2(\Pi)$ of dimension $2p+3$.

\begin{rem}
{\rm Our first approach to these problems was to use the spectral sequence associated to our
complex, filtred by the valuation (whith respect to the
quasihomogeneous degree). 
But the 
method we present here gives better results.}
\end{rem}
\endrem

\end{document}